\documentclass[oneside,english]{amsart}

\usepackage{amsmath}
\usepackage{amsthm}
\usepackage{amssymb}
\usepackage{amscd}
\usepackage{xspace}
\usepackage{verbatim}
\usepackage{geometry}
\usepackage{esint}
\usepackage{geometry}
\usepackage{fancyhdr}
\usepackage{xargs}
\usepackage{xifthen}
\usepackage{hyperref}

\newcommandx{\Set}[2][2=]{
    \ifthenelse{\isempty{#2}}
        {\left\{ {#1} \right\}}
        {\left\{  {#1}  \, \middle| \, {#2} \right\}}
}

\newtheorem{theorem}{Theorem}[section]
\newtheorem*{theorem*}{Theorem}

\newtheorem{lemma}{Lemma}[section]
\newtheorem{proposition}{Proposition}[section]
\newtheorem*{proposition*}{Proposition}

\theoremstyle{definition}
\newtheorem{definition}{Definition}[section]
\newtheorem*{acknowledgment}{Acknowledgment}
\theoremstyle{remark}
\newtheorem{remark}{Remark}[section]
\numberwithin{equation}{section}

\newcommand{\e}{\varepsilon}

\newcommand{\N}{\mathbb N}

\newcommand{\R}{\mathbb{R}}

\DeclareMathOperator{\supp}{supp}

\DeclareMathOperator{\cp}{Cap}

\DeclareMathOperator{\loc}{loc}

\newcommandx{\norm}[1][1=\cdot]{\left|{#1}\right|}
\newcommandx{\supnorm}[2][1=\cdot,2=]{\nnorm[#1]_{\infty,#2}}
\newcommandx{\nnorm}[1][1=\cdot]{\left|\left|{#1}\right|\right|}

\date{\today}

\begin{document}

\title[The refined area formula for Sobolev mappings $W^{k,p}$]{The Refined Area Formula for Sobolev Mappings $W^{k,p}$}

\author{Paz Hashash}

\maketitle
\begin{center}
\textit{Dedicated to Vladimir Gol'dstein on the occasion of his 80th birthday.}
\end{center}

\begin{abstract}
We establish the area formula for change-of-variable mappings in the Sobolev space $W^{k,p}_{\text{loc}}$. Our approach relies on constructing Lipschitz approximations of Sobolev functions that agree with the original functions outside a set of Riesz capacity zero.
\end{abstract}

\section{\textbf{Introduction}}

In this paper, we study approximations of Sobolev functions in the space $W^{k,p}_{\text{loc}}(\Omega)$ by Lipschitz continuous functions, where $k > 1$ is a natural number, $1 \leq p < \infty$, and $\Omega \subset \mathbb{R}^n$ is an open set. These approximations have applications in geometric measure theory~\cite{Hajlasz1993,Vodop'yanovUkhlov1996} and in the analysis of Sobolev functions on metric measure spaces~\cite{HajlaszKinnunen1998,Heinonen2001,HeinonenKoskela2015}. 

Our primary motivation is to prove the area formula for Sobolev functions in $W^{k,p}_{\text{loc}}(\Omega)$, where the change-of-variable mapping is defined almost everywhere with respect to Riesz capacity $R_{k-1,p}$.

The area formula for Lipschitz mappings was first proved by Federer; see~\cite{Fe69} or \cite{AFP2000}. It states the following: let $A \subset \mathbb{R}^n$ be a Lebesgue measurable set, and let $\varphi : A \to \mathbb{R}^n$ be a Lipschitz mapping. Then, for any non-negative Lebesgue measurable function $f : A \to \mathbb{R}$, we have
\begin{equation}
\label{eq: area formula for Lipschitz mappings}
\int_{A} f(x) |J(x, \varphi)| \, dx = \int_{\mathbb{R}^n} \left( \sum_{x \in \varphi^{-1}\{y\} \cap A} f(x) \right) \, dy,
\end{equation}
where $J(x, \varphi) = \det(D\varphi(x))$ denotes the Jacobian of $\varphi$ at $x$, and $D\varphi(x)$ is interpreted as the approximate differential. In particular, the function 
$y\mapsto  \sum_{x \in \varphi^{-1}\{y\} \cap A} f(x)$ is Lebesgue measurable on $A$.

We prove a generalized version of this formula under weaker assumptions:

\begin{theorem}
\label{thm:general area formula}
Let $A \subset \mathbb{R}^n$ be a Lebesgue measurable set and let $\varphi : A \to \mathbb{R}^n$ be a Lebesgue measurable mapping. Suppose there exists a sequence of Lebesgue measurable sets $\{A_m\}_{m=1}^\infty \subset A$, with $A_m \subset A_{m+1}$, and $\varphi|_{A_m}$ is Lipschitz continuous for every $m \in \mathbb{N}$. We denote $S := A \setminus \bigcup_{m=1}^\infty A_m$. Then, for any non-negative Lebesgue measurable function $f : A \to \mathbb{R}$, we have:
\begin{equation}
\label{eq:formula 3}
\int_{A\setminus S} f(x) |J(x, \varphi)| \, dx = \int_{\mathbb{R}^n} \left( \sum_{x \in \varphi^{-1}\{y\} \cap (A \setminus S)} f(x) \right) \, dy.
\end{equation}
\end{theorem}

Formula \eqref{eq:formula 3} follows from the classical area formula for Lipschitz mappings and the monotone convergence theorem. Using Theorem \ref{thm:general area formula}, we prove the following main result:

\begin{theorem}
\label{thm:area formula for Riesz capacity}
Let $k > 1$ be a natural number and let $1 \leq p < \infty$ be such that $(k-1)p < n$. Let $\Omega \subset \mathbb{R}^n$ be an open set and let $\varphi \in W^{k,p}_{\text{loc}}(\Omega, \mathbb{R}^n)$. Then there exists a set $S \subset \Omega$ such that $R_{k-1,p}(S) = 0$, and for any non-negative Lebesgue measurable function $f : \Omega \to \mathbb{R}$, we have:
\begin{equation}
\int_{\Omega} f(x) |J(x, \varphi)| \, dx = \int_{\mathbb{R}^n} \left( \sum_{x \in \varphi^{-1}\{y\} \cap (\Omega \setminus S)} f(x) \right) \, dy.
\end{equation}
\end{theorem}

An area formula theorem was previously proved by Reshetnyak, Vodop'yanov, and Gol'dstein for Sobolev mappings in $W^{1,n}$ with respect to the Lebesgue measure in~\cite{goldstein1979geometric}. In the case of $W^{1,1}_{\text{loc}}$ with the Lebesgue measure, the area formula was proved by Piotr Hajłasz in~\cite{Hajlasz1993}; see also~\cite{Hajlasz1996}. Theorem \ref{thm:area formula for Riesz capacity} was established for $k = 2$ and the Sobolev capacity $\cp_p$ in~\cite{HashashUkhlov2022}. 

To prove Theorem~\ref{thm:area formula for Riesz capacity}, we establish the following Lipschitz approximations for mappings in $W^{k,p}_{\text{loc}}$ involving Riesz capacity:

\begin{theorem}
\label{thm:LipCap}
Let $\Omega \subset \mathbb{R}^n$ be an open set, let $k \in \mathbb{N}$ with $k > 1$, and let $1 \leq p < \infty$ be such that $(k-1)p < n$. Suppose $f \in W^{k,p}_{\loc}(\Omega)$. Then there exists a sequence of closed sets $\{C_l\}_{l=1}^{\infty}$ such that for every $l = 1, 2, \dots$, we have $C_l \subset C_{l+1} \subset \Omega$, and the restriction $f^* \vert_{C_l}$ is a Lipschitz continuous function defined $R_{k-1,p}$-almost everywhere in $C_l$, and
\begin{equation}
R_{k-1,p}\left(\Omega\setminus\bigcup_{l=1}^{\infty}C_l\right)=0.
\end{equation}
Here $f^*$ is the precise representative of the function $f$, see Definition \ref{def:precise representative}. 
\end{theorem}
See also a related approximation result in \cite{Maly1993}.

We would like to provide some motivation for Theorems~\ref{thm:area formula for Riesz capacity} and~\ref{thm:LipCap}. As previously mentioned, the area formula for Sobolev functions in \( W^{1,1}_{\text{loc}} \) was established by Piotr Hajłasz in~\cite{Hajlasz1993}. To derive the area formula, he showed that such Sobolev functions can be approximated by Lipschitz continuous functions outside a set of Lebesgue measure zero. However, sets of Lebesgue measure zero can still be geometrically large. For instance, in \( \mathbb{R}^3 \), a two-dimensional surface has Lebesgue measure zero, so this approximation result does not directly allow one to restrict a Sobolev function to a two-dimensional surface (when needed).

More precisely, let us consider a mapping \( \varphi \in W^{2,2}(\mathbb{R}^3, \mathbb{R}^3) \). Theorem~\ref{thm:LipCap} implies that \( \varphi \) is defined outside a set \( S \subset \mathbb{R}^3 \) of \( R_{1,2} \)-capacity zero. A classical result (see Theorem~2.6.16 in~\cite{Ziemer2012}) asserts that \( \mathcal{H}^2(S) = 0 \), where \( \mathcal{H}^2 \) denotes the two-dimensional Hausdorff measure. Hence, \( \varphi \) can be meaningfully restricted to any two-dimensional surface in \( \mathbb{R}^3 \), since such surfaces intersect \( S \) in a set of \( \mathcal{H}^2 \)-measure zero.

Theorem~\ref{thm:area formula for Riesz capacity} strengthens Hajłasz’s area formula in~\cite{Hajlasz1993} by showing that the change-of-variable mapping \( \varphi \) is defined almost everywhere not only with respect to the Lebesgue measure, but also with respect to the Riesz capacity. Recall that Lebesgue measure is absolutely continuous with respect to the Riesz capacity (see Remark \ref{rem:Riesz capacity is an outer measure}). In particular, \( \varphi \) is defined on a domain that is larger in the capacitary sense.
\\

The article is organized as follows: In Section~\ref{sec:the arear formula}, we prove the area formula in a general setting under the assumption of the existence of Lipschitz approximations. In Section~\ref{sec:Lebesgue points of Soblev functions and capacities}, we recall results on Lebesgue points of Sobolev functions. Finally, in Section~\ref{sec:Lipschitz approximations}, we prove the required Lipschitz approximation for Sobolev functions that yields Theorem~\ref{thm:area formula for Riesz capacity}.

\section{\textbf{The area formula}}
\label{sec:the arear formula}

\begin{proof}[\textbf{Proof of Theorem \ref{thm:general area formula}}]
By the area formula for the Lipschitz mapping $\varphi$ on the set $A_m$, $m \in \mathbb{N}$ (see equation~\eqref{eq: area formula for Lipschitz mappings}), we obtain, for any non-negative Lebesgue measurable function $f:A \to \mathbb{R}$,
\begin{equation}
\label{eq:formula 2'}
\int_{A_m} f(x) |J(x,\varphi)|\,dx = \int_{\mathbb{R}^n} \left( \sum_{x \in \varphi^{-1}\{y\} \cap A_m} f(x) \right)\,dy.
\end{equation} 
Observe that the sequence of non-negative functions $f\chi_{A_m}$ is non-decreasing on $A$, and converges almost everywhere (with respect to Lebesgue measure) to $f\chi_{A \setminus S}$. Here, $\chi_E$ denotes the characteristic function of the set $E \subset \mathbb{R}^n$. 

Note also that the functions 
\[
g_m(y) := \sum_{x \in \varphi^{-1}\{y\} \cap A_m} f(x)
\]
form a monotone non-decreasing sequence of non-negative Lebesgue measurable functions. Their measurability follows from the area formula.

Therefore, applying the monotone convergence theorem to both sides of equation~\eqref{eq:formula 2'}, we obtain
\begin{multline}
\label{eq:formula 4}
\int_{A \setminus S} f(x) |J(x,\varphi)|\,dx 
= \int_{A} f(x)\chi_{A \setminus S}(x)|J(x,\varphi)|\,dx 
= \int_{A} \left( \lim_{m \to \infty} f(x) \chi_{A_m}(x) \right)|J(x,\varphi)|\,dx \\
= \lim_{m \to \infty} \int_{A} f(x) \chi_{A_m}(x)|J(x,\varphi)|\,dx 
= \lim_{m \to \infty} \int_{A_m} f(x)|J(x,\varphi)|\,dx 
\\
= \lim_{m \to \infty} \int_{\mathbb{R}^n} \left( \sum_{x \in \varphi^{-1}\{y\} \cap A_m} f(x) \right)\,dy 
= \int_{\mathbb{R}^n} \lim_{m \to \infty} \left( \sum_{x \in \varphi^{-1}\{y\} \cap A_m} f(x) \right)\,dy 
\\
\int_{\mathbb{R}^n} \lim_{m \to \infty} \left( \sum_{x \in \varphi^{-1}\{y\}} f(x)\chi_{A_m}(x) \right)\,dy 
= \int_{\mathbb{R}^n} \left( \sum_{x \in \varphi^{-1}\{y\}} f(x)\chi_{A \setminus S}(x) \right)\,dy 
= \int_{\mathbb{R}^n} \left( \sum_{x \in \varphi^{-1}\{y\} \cap (A \setminus S)} f(x) \right)\,dy.
\end{multline}
\end{proof}

\section{\textbf{Lebesgue points of Sobolev functions and capacities}}
\label{sec:Lebesgue points of Soblev functions and capacities}

\begin{definition}[Sobolev Space \( W^{k,p}(\Omega) \)]
Let \( \Omega \subset \mathbb{R}^n \) be an open set, \( k \in \mathbb{N} \), and \( 1 \leq p \leq \infty \). The Sobolev space \( W^{k,p}(\Omega) \) is defined by
\[
W^{k,p}(\Omega) := \Set{ f \in L^p(\Omega)}[D^\alpha f \in L^p(\Omega) \text{ for all multi-indices } \alpha \text{ with } |\alpha| \leq k ],
\]
where \( D^\alpha f \) denotes the weak derivative of \( f \) of order \( \alpha \), and \( |\alpha| \) is the order of the multi-index \( \alpha \).

The Sobolev norm on \( W^{k,p}(\Omega) \) is given by
\[
\|f\|_{W^{k,p}(\Omega)} := \left( \sum_{|\alpha| \leq k} \|D^\alpha f\|_{L^p(\Omega)}^p \right)^{1/p}, \quad \text{for } 1 \leq p < \infty,
\]
and
\[
\|f\|_{W^{k,\infty}(\Omega)} := \max_{|\alpha| \leq k} \|D^\alpha f\|_{L^\infty(\Omega)}.
\]

The local Sobolev space \( W^{k,p}_{\text{loc}}(\Omega) \) is defined as the space of all functions \( f \in L^1_{\text{loc}}(\Omega) \) such that for every open and bounded set \( U \subset \Omega \) whose closure satisfies \( \overline{U} \subset \Omega \), we have \( f \in W^{k,p}(U) \).

We denote by \( \nabla f \) the weak gradient of a function \( f \in W^{1,p}(\Omega) \).
\end{definition}

For more information about Sobolev spaces see, for example, \cite{EvansGariepy1992,Mazya2011,Ziemer2012}.

The following two Lemmas are taken from Ziemer's book \cite{Ziemer2012}. We will use them in the proof of Theorem \ref{thm:Lebesgue points of Sobolev functins}. 
\begin{lemma}[Lemma 3.1.1 in \cite{Ziemer2012}]
\label{lem:Lemma 3.1.1}
Let $x_0 \in \mathbb{R}^n$, $r\in (0,\infty)$, and $B(x_0,r)\subset \mathbb{R}^n$ be a ball. Let $1\leq p\leq \infty$ and $f \in W^{1,p}(B(x_0,r))$. Let $0 < \delta < r$. Then
\begin{multline}
r^{-n} \int_{B(x_0,r)} f(y) \, dy - \delta^{-n} \int_{B(x_0,\delta)} f(y) \, dy 
= \frac{1}{n} r^{-n} \int_{B(x_0,r)} [\nabla f(y) \cdot (y - x_0)] \, dy 
\\
\quad - \frac{1}{n} \delta^{-n} \int_{B(x_0,\delta)} [\nabla f(y) \cdot (y - x_0)] \, dy 
\\
\quad - \frac{1}{n} \int_{B(x_0,r) \setminus B(x_0,\delta)} |y - x_0|^{-n} [\nabla f(y) \cdot (y - x_0)] \, dy.
\end{multline}
\end{lemma}

\begin{lemma}[Lemma 3.1.3 in \cite{Ziemer2012}]
\label{lem:Lemma 3.1.3}
Let $\ell$ be a positive real number, $k$ a positive integer and $p\geq 1$ such that $(k + \ell - 1)p < n$. Then there exists a constant $C = C(n, k, \ell)$ such that
\begin{equation}
\int_{\mathbb{R}^n} |y - x|^{\ell - n} |f(y)| \, dy 
\leq C \sum_{|\alpha| = k} \int_{\mathbb{R}^n} |y - x|^{\ell - n + k} |D^\alpha f(y)| \, dy
\end{equation}
for all $x\in \mathbb{R}^n$ and $f\in W^{k,p}(\mathbb{R}^n)$.
\end{lemma}

Recall the definition of Riesz capacity:
\begin{definition}
Let $0 < \alpha < n$. The {\it Riesz kernel} is defined as the function
\begin{equation}
I_\alpha(x) := \frac{1}{\gamma(\alpha)}\frac{1}{|x|^{n-\alpha}},
\end{equation}
where $\gamma(\alpha) = \pi^{n/2} 2^\alpha \frac{\Gamma(\alpha/2)}{\Gamma(n/2 - \alpha/2)}$, and $\Gamma(t):=\int_0^\infty e^{-x}x^{t-1}dx$, $t\in (0,\infty)$, denotes the gamma function.
Let $\varphi:\mathbb{R}^n\to \mathbb{R}$ be a non-negative Lebesgue measurable function.
The {\it Riesz potential} of $\varphi$ is given by the convolution $I_\alpha * \varphi(x)$, which is defined as
\begin{equation}
I_\alpha * \varphi(x) = \frac{1}{\gamma(\alpha)} \int_{\mathbb{R}^n} \frac{\varphi(y)}{|x - y|^{n-\alpha}} \, dy.
\end{equation} 
For $\alpha \in (0,\infty)$ and $p \in [1,\infty)$ with $\alpha p < n$, the {\it Riesz capacity} $R_{\alpha,p}$ is defined for any subset $E \subset \mathbb{R}^n$ as
\[
R_{\alpha,p}(E) = \inf \left\{ \|\varphi\|_{L^p(\mathbb{R}^n)}^p : I_\alpha * \varphi(x) \geq 1 \text{ on } E, \, \varphi \in L^p(\mathbb{R}^n), \, \varphi \geq 0 \right\}.
\]
\end{definition}

\begin{remark}
\label{rem:Riesz capacity is an outer measure}
It is known that the Riesz capacity is an outer measure and that the Lebesgue measure is absolutely continuous with respect to the Riesz capacity $R_{\alpha,p}$. As the author could not find a reference establishing these properties for the full range $1 \leq p < \infty$, we provide an explanation here. Hereafter, we denote the $n$-dimensional Lebesgue measure of a set $E \subset \mathbb{R}^n$ by $|E|$.

Let $1 \leq p < \infty$ and $\alpha \in (0,\infty)$ such that $\alpha p < n$, and set $p^* := \frac{np}{n - \alpha p}$.

The Lebesgue measure is absolutely continuous with respect to the Riesz capacity $R_{\alpha,p}$. Indeed, let $f \in L^p(\mathbb{R}^n)$, $f \geq 0$, and suppose that $I_\alpha * f \geq 1$ on a set $E$. In case $1 < p < \infty$
\begin{equation}
\label{eq:estimate of Lebesgue measure in case p>1}
|E|^{\frac{p}{p^*}} = \left(\int_E 1\, dx\right)^{\frac{p}{p^*}} 
\leq \left( \int_E (I_\alpha * f(x))^{p^*} dx \right)^{\frac{p}{p^*}} 
\leq \| I_\alpha * f \|_{L^{p^*}(\mathbb{R}^n)}^p 
\leq C(n,p) \|f\|_{L^p(\mathbb{R}^n)}^p.
\end{equation}
For the last inequality in \eqref{eq:estimate of Lebesgue measure in case p>1}, see Theorem 2.8.4 in \cite{Ziemer2012}. This inequality is proved using the Hardy–Littlewood inequality, which is valid for $1 < p \leq \infty$, but fails for $p = 1$; see Theorem 2.8.2 in \cite{Ziemer2012} for a proof of this inequality.
By the definition of $R_{\alpha,p}$, inequality \eqref{eq:estimate of Lebesgue measure in case p>1} shows that if $R_{\alpha,p}(E) = 0$, then $|E| = 0$.

For the case $p = 1$, recall first that the Riesz kernel $I_\alpha$ lies in the Lorentz space (or weak Lebesgue space) $L^{\frac{n}{n - \alpha}, \infty}(\mathbb{R}^n)$. For a discussion of this fact, see the end of Section 2.10.1 in \cite{Ziemer2012}. Recall that for a Lebesgue measurable function $g:\R^n\to \R$ and positive number $\gamma$
\begin{equation} 
\|g\|^{\gamma}_{L^{\gamma, \infty}(\mathbb{R}^n)}:=\sup\limits_{0<t<\infty}t\left| \left\{ x \in \mathbb{R}^n : |g(x)|^\gamma \geq t \right\} \right|.
\end{equation}
It follows that
\begin{equation}
\label{eq:estimate of Lebesgue measure in case p=1}
|E|^{\frac{n - \alpha}{n}} 
\leq \left| \left\{ x \in \mathbb{R}^n : (I_\alpha * f(x))^{\frac{n}{n - \alpha}} \geq 1 \right\} \right|^{\frac{n - \alpha}{n}} 
\leq \| I_\alpha * f \|_{L^{\frac{n}{n - \alpha}, \infty}(\mathbb{R}^n)} 
\leq C(n,\alpha) \|I_\alpha\|_{L^{\frac{n}{n - \alpha}, \infty}(\R^n)} \|f\|_{L^1(\mathbb{R}^n)}.
\end{equation}
The last inequality in \eqref{eq:estimate of Lebesgue measure in case p=1} follows from Young's inequality for weak-type spaces; see Theorem 1.2.13 in \cite{GrafakosClassical}. Hence, if $R_{\alpha,1}(E) = 0$, then $|E| = 0$.

The fact that the Riesz capacity is an outer measure follows from its definition. Suppose $E \subset \bigcup_{i=1}^\infty E_i$. We aim to prove that
\[
R_{\alpha,p}(E) \leq \sum_{i=1}^\infty R_{\alpha,p}(E_i).
\]
Assume that $\sum_{i=1}^\infty R_{\alpha,p}(E_i) < \infty$. For each $i \in \mathbb{N}$, choose a function $f_i \in L^p(\mathbb{R}^n)$ with $f_i \geq 0$, $I_\alpha * f_i \geq 1$ on $E_i$, and $\|f_i\|_{L^p}^p \leq R_{\alpha,p}(E_i) + \frac{\epsilon}{2^i}$. Define $f := \sup_{i \in \mathbb{N}} f_i$. Then $f \in L^p(\mathbb{R}^n)$, $f \geq 0$, $I_\alpha * f \geq 1$ on $E$, and
\[
R_{\alpha,p}(E) \leq \|f\|_{L^p(\mathbb{R}^n)}^p \leq  \sum_{i=1}^\infty \|f_i\|_{L^p(\mathbb{R}^n)}^p \leq \sum_{i=1}^\infty \left(R_{\alpha,p}(E_i) + \frac{\epsilon}{2^i} \right) = \sum_{i=1}^\infty R_{\alpha,p}(E_i) + \epsilon.
\]
Since $\epsilon > 0$ is arbitrary, the result follows.
\end{remark}

The following theorem is attributed to Ziemer \cite{Ziemer2012} (see Theorem 3.1.4.), where it is assumed that $p>1$. We present Ziemer's proof here to demonstrate that the same argument remains valid in the case $p = 1$.
\begin{theorem}
\label{thm:Lebesgue points of Sobolev functins}
Let $1\leq p<\infty$, $k$ be a positive integer such that $kp < n$, and let $f\in W^{k,p}(\mathbb{R}^n)$. Define
\begin{equation}
\label{eq:definition of the non Lebesgue points of Ries potential of high derivatives}
E:=\Set{x\in\mathbb{R}^n}[I_k*g(x)=\infty],\quad g(y) := \sum_{|\alpha| = k} |D^\alpha f(y)|.
\end{equation}
Then $R_{k,p}(E) = 0$ and the limit $\lim_{\delta \to 0^+} \fint_{B(x,\delta)} f(y) \, dy$ exists and belongs to $\mathbb{R}$ for every $x\in \mathbb{R}^n\setminus E$.
\end{theorem}

\begin{proof}
We prove first that $R_{k,p}(E) = 0$. Let $\epsilon\in (0,\infty)$. Then, the function $\epsilon g$ belongs to $L^p(\mathbb{R}^n)$ and for every $x\in E$ we get $I_k*(\epsilon g)(x)=\infty>1$. Therefore, by the definition of $R_{k,p}$, we get $R_{k,p}(E)\leq \epsilon^p\|g\|^p_{L^p(\mathbb{R}^n)}$. Since $\epsilon$ is arbitrary, we conclude that $R_{k,p}(E) = 0$.

Consider $x \in \mathbb{R}^n \setminus E$. By Lemma \ref{lem:Lemma 3.1.1} with $r=1$ and $0<\delta<1$, we get 
\begin{multline}
\label{eq:usage of lemma 4.1}
\int_{B(x,1)} f(y) \, dy - \delta^{-n} \int_{B(x,\delta)} f(y) \, dy
= \frac{1}{n} \int_{B(x,1)} [\nabla f(y) \cdot (y - x)] \, dy
\\
- \frac{1}{n} \delta^{-n} \int_{B(x,\delta)} [\nabla f(y) \cdot (y - x)] \, dy
- \frac{1}{n} \int_{B(x,1)\setminus B(x,\delta)} |y - x|^{-n} [\nabla f(y) \cdot (y - x)] \, dy.
\end{multline}
We will prove that the limit on the right-hand side of the equation in \eqref{eq:usage of lemma 4.1} exists as $\delta \to 0^+$ and, consequently, deduce the existence of the limit on the left-hand side of the equation in \eqref{eq:usage of lemma 4.1}.

We first prove that: 
\begin{equation}
\label{eq:finiteness of derivative on a ball}
\int_{B(x,1)} |y - x|^{1-n} |\nabla f(y)| \, dy < \infty.
\end{equation}
When $k = 1$ ($k$ is the differentiability degree of $f$), \eqref{eq:finiteness of derivative on a ball} follows from \eqref{eq:definition of the non Lebesgue points of Ries potential of high derivatives} because for $x\notin E$, we get
\begin{equation}
\infty>I_1*g(x)= \frac{1}{\gamma(1)}\sum_{|\alpha| = 1} \int_{\mathbb{R}^n}|x-y|^{1-n}|D^\alpha f(y)|dy\geq \frac{1}{\gamma(1)}\int_{B(x,1)} |y - x|^{1-n} |\nabla f(y)| \, dy.
\end{equation}
When $k > 1$, it follows from Lemma \ref{lem:Lemma 3.1.3} with $\ell = 1$ and $k - 1$ substituted for $k$, and multi-index $\beta$ such that $|\beta|=1$, 
\begin{equation}
\label{eq:Riesz potential is dominated by Riesz potential of higher derivaties}
\int_{\mathbb{R}^n} |y - x|^{1-n} |D^\beta f(y)| \, dy \leq C \sum_{|\alpha| = k-1} \int_{\mathbb{R}^n} |y - x|^{k-n} |D^\alpha \left(D^\beta f(y)\right)| \, dy.
\end{equation}
Thus, taking the sum over all such multi-indices $\beta$ in \eqref{eq:Riesz potential is dominated by Riesz potential of higher derivaties}, we obtain
\begin{equation}
\int_{\mathbb{R}^n} |y - x|^{1-n} |\nabla f(y)| \, dy\leq C\,\sum_{|\alpha| = k} \int_{\mathbb{R}^n} |y - x|^{k-n} |D^\alpha f(y)| \, dy=C\gamma(k)I_k*g(x)<\infty.
\end{equation}
Thus \eqref{eq:finiteness of derivative on a ball} still holds when $k > 1$. This proves \eqref{eq:finiteness of derivative on a ball}.  

It follows from \eqref{eq:finiteness of derivative on a ball} that
\begin{equation}
\label{eq:limit of du times |y-x|}
\lim_{\delta \to 0^+} \int_{B(x,\delta)} |y - x|^{1-n} |\nabla f(y)| \, dy =\int_{\{x\}} |y - x|^{1-n} |\nabla f(y)| \, dy= 0.
\end{equation}
We prove now that the following limit exists:
\begin{equation}
\label{eq:usage of lemma 4.1(1)}
\lim_{\delta \to 0^+} \int_{B(x,1)\setminus B(x,\delta)} |y - x|^{-n} [\nabla f(y) \cdot (y - x)] \, dy.
\end{equation}
Indeed, since
\begin{multline}
\label{eq:estimate for proving limit}
\left|\int_{B(x,1)\setminus B(x,\delta)} |y - x|^{-n} [\nabla f(y) \cdot (y - x)] \, dy
-\int_{B(x,1)} |y - x|^{-n} [\nabla f(y) \cdot (y - x)] \, dy\right|
\\
=\left|\int_{B(x,\delta)} |y - x|^{-n} [\nabla f(y) \cdot (y - x)] \, dy\right|
\leq \int_{B(x,\delta)} |y - x|^{1-n} |\nabla f(y)|dy,
\end{multline}
the existence of the limit in \eqref{eq:usage of lemma 4.1(1)} follows from \eqref{eq:estimate for proving limit} and \eqref{eq:limit of du times |y-x|}. 

From \eqref{eq:limit of du times |y-x|}, we obtain 
\begin{equation}
\label{eq:usage of lemma 4.1 (2)}
\lim_{\delta \to 0^+}\left| \delta^{-n} \int_{B(x,\delta)} [\nabla f(y) \cdot (y - x)] \, dy\right|\leq \lim_{\delta \to 0^+} \int_{B(x,\delta)} |y - x|^{1-n} |\nabla f(y)| \, dy=0.
\end{equation}
It now follows from \eqref{eq:usage of lemma 4.1}, \eqref{eq:usage of lemma 4.1(1)}, and \eqref{eq:usage of lemma 4.1 (2)} that the limit $\lim_{\delta\to 0^+}\delta^{-n} \int_{B(x,\delta)} f(y) \, dy$ exists and belongs to $\R$.

\end{proof}

The following theorem was proved in \cite{Tai2001}:
\begin{theorem}
\label{thm:weak type inequality for capacities}
Let \( k \) be a positive integer with \( kp < n \), where \( 1 \leq p < \infty \). Then, there exists a positive constant \( C \) depending on \( n \), \( k \), and \( p \) such that 
\[
R_{k,p}\left(\Set{x \in \mathbb{R}^n}[Mf(x) > \lambda]\right) \leq \frac{C}{\lambda^p}\|f\|_{W^{k,p}(\mathbb R^n)}^p
\]
for \( f \in W^{k,p}(\mathbb{R}^n) \) and \( \lambda\in (0,\infty)\), where \( R_{k,p} \) is the Riesz capacity and
\[
Mf(x) = \sup_{r\in (0,\infty)} \frac{1}{|B(x,r)|} \int_{B(x,r)} |f(y)| \, dy.
\]
\end{theorem}

\begin{remark}
\label{rem:chebyshev inequality for Riesz capacity is valid for vector valued Sobolev mappings}
Note that Theorem \ref{thm:weak type inequality for capacities} is also valid for vector-valued Sobolev mappings: let $F:\R^n\to \R^m$, $F=(F_1,\dots,F_m)$ be a Lebesgue measurable mapping. Since $|F(y)|\leq \sum_{i=1}^m|F_i(y)|$ for a.e $y\in \R^n$, we get
\begin{equation}
M|F|(x)=\sup_{r\in (0,\infty)}  \fint_{B(x,r)} |F(y)|\,dy
\leq \sum_{i=1}^m\sup_{r\in (0,\infty)}  \fint_{B(x,r)} |F_i(y)|\,dy
= \sum_{i=1}^m MF_i(x).
\end{equation}
Therefore, for $\lambda\in (0,\infty)$,
\begin{equation}
\left\{M|F|>\lambda\right\}\subset \left\{\sum_{i=1}^m MF_i>\lambda\right\}\subset \bigcup_{i=1}^m\left\{ MF_i>\frac{\lambda}{m}\right\}.
\end{equation}
Thus, if $F\in W^{k,p}(\R^n,\R^m)$, then we get from Theorem \ref{thm:weak type inequality for capacities}
\begin{multline}
R_{k,p}\left(\Set{x \in \mathbb{R}^n}[ M|F|(x) > \lambda]\right)
\leq \sum_{i=1}^m R_{k,p}\left(\Set{x \in \mathbb{R}^n}[MF_i(x) > \frac{\lambda}{m}]\right)
\\
\leq \frac{C}{\lambda^p}\sum_{i=1}^m\|F_i\|_{W^{k,p}(\mathbb{R}^n)}^p
\leq \frac{C}{\lambda^p}\|F\|_{W^{k,p}(\mathbb{R}^n,\R^m)}^p,
\end{multline}
where $C$ depends only on $k,p,n,m$.
\end{remark}

\section{\textbf{Lipschitz approximations of Sobolev functions quasi-everywhere}}
\label{sec:Lipschitz approximations}
\begin{definition}
\label{def:precise representative}
Let $\Omega\subset\mathbb R^n$ be an open set and  
$f\in L^1_{\loc}(\Omega).$
The {\it precise representative} of 
$f$ is defined by 
\begin{equation}
\label{precise representative}
f^*:\Omega\to \mathbb R,\quad f^*(x):=
\begin{cases}
\lim_{\epsilon\downarrow 0}f_{B(x,\epsilon)},\quad & \text{if the limit exists and belongs to $\mathbb R$};\\
0, \quad & \text{otherwise};
\end{cases}
\end{equation}
where
\begin{equation}
f_{B(x,\epsilon)}:=\fintop_{B(x,\epsilon)}f(y)dy=\frac{1}{\left|B(x,\epsilon)\right|}\intop_{B(x,\epsilon)}f(y)dy.
\end{equation}
\end{definition}

\begin{lemma}
\label{lem:level sets of the H.L function are open}
Let $f\in L^1_{\loc}(\mathbb R^n)$. Let 
\begin{equation}
\label{eq:H.L maximal function 1}
Mf(x):=\sup_{r\in(0,\infty)}\fint_{B(x,r)}|f(z)|dz=\sup_{r\in(0,\infty)}|f|_{B(x,r)}
\end{equation}
be the Hardy-Littlewood maximal function, and let us denote for every $\lambda\in (0,\infty)$
\begin{equation}
A_\lambda:=\Set{x\in \mathbb R^n}[Mf(x)\leq \lambda].
\end{equation}
Then $A_\lambda$ is a closed set.
\end{lemma}

\begin{proof}
Let $\{x_i\}_{i=1}^\infty \subset A_\lambda$ be a sequence of points such that $x_i\to x$ as $i\to \infty$. We prove that $x\in A_\lambda$, and, therefore, $A_\alpha$ is closed. For every $r>0$, the function $z\mapsto |f|_{B(z,r)}$ is continuous at $x$. Therefore, for $\varepsilon>0$, there exists $\delta>0$ such that $\left||f|_{B(z,r)}-|f|_{B(x,r)}\right|\leq \varepsilon$ for all $z\in B(x,\delta)$. 
Thus, for sufficiently large indices $i$ such that $x_i\in B(x,\delta)$, we have $|f|_{B(x,r)}\leq \varepsilon+|f|_{B(x_i,r)}\leq \varepsilon+\lambda$. Consequently, $\sup_{r\in (0,\infty)}|f|_{B(x,r)}\leq \varepsilon+\lambda$. Since $\varepsilon>0$ is arbitrary, we conclude that $x\in A_\lambda$, proving that $
A_\lambda$ is a closed set.
\end{proof}

\begin{theorem}[The local Poincar\'e inequality for functions of $W^{1,1}_{\loc}(\mathbb R^n)$]
\label{thm:Poincare inequality}
There exists a constant 
$C=C(n)$ 
such that 
\begin{equation}
\fintop_{B(x,r)}|f(y)-f_{B(x,r)}|dy\leq Cr\fintop_{B(x,r)}|\nabla f(y)|dy,
\end{equation}
for every ball $B(x,r)\subset \mathbb R^n$ and every 
$f\in W^{1,1}(B(x,r))$.       
\end{theorem} 
A proof of Theorem \ref{thm:Poincare inequality} can be found in many books; see, for example, Theorem 2 in Section 4.5.2 of \cite{EvansGariepy1992}.

\begin{lemma}
\label{lem:Lipschitz approximations}
Let $1\leq p\leq \infty$, and $f\in W^{1,p}(\mathbb{R}^n)$. For every $\alpha\in (0,\infty)$, 
define 
\begin{equation}
A_\alpha = \Set{x\in \mathbb{R}^n \mid M|\nabla f|(x)\leq \alpha}.
\end{equation}
Let $\mathcal{N}\subset \mathbb{R}^n$ denote the set of points $x\in \mathbb{R}^n$ such that the limit $\lim_{r\to 0^+}f_{B(x,r)}$ does not exist or does not belong to $\mathbb{R}$.  
Then, the restriction 
$f^* \vert_{A_\alpha \setminus \mathcal{N}}$ 
is a Lipschitz continuous function on $A_\alpha \setminus \mathcal{N}$. 
\end{lemma}

\begin{proof}
We begin by proving that points \(x \in A_\alpha\) satisfy the following average approximation property: if \(x \in A_\alpha\), then, by the Poincar\'e inequality (Theorem~\ref{thm:Poincare inequality}), for all positive numbers \(r > s\), we have
\begin{equation}
\label{eq:estimate for differences between averages, in case W11(1)}
|f_{B(x,r)} - f_{B(x,s)}|
\leq \left(\frac{r}{s}\right)^n C(n) r \alpha,
\end{equation}
where \(C(n)\) is a constant depending only on the dimension \(n\). Indeed,
\begin{multline}
\label{eq:estimate for differences between averages, in case W11(2)}
|f_{B(x,r)} - f_{B(x,s)}|
\leq \fintop_{B(x,s)} |f(y) - f_{B(x,r)}| \, dy
= \frac{|B(x,r)|}{|B(x,s)|} \cdot \frac{1}{|B(x,r)|} \intop_{B(x,s)} |f(y) - f_{B(x,r)}| \, dy \\
\leq \left(\frac{r}{s}\right)^n \fintop_{B(x,r)} |f(y) - f_{B(x,r)}| \, dy
\leq \left(\frac{r}{s}\right)^n C(n) r \fintop_{B(x,r)} |\nabla f(y)| \, dy
\leq \left(\frac{r}{s}\right)^n C(n) r \alpha.
\end{multline}

Next, we use the average approximation~\eqref{eq:estimate for differences between averages, in case W11(1)} to show that one can approximate the precise representative of the function \(f\) at the point \(x \in A_\alpha\) by an average of \(f\) around \(x\). More precisely, we prove that: if \(x \in A_\alpha \setminus \mathcal{N}\), then
\begin{equation}
\label{eq:estimate for |f^*(x)-f_B(x,r)|,in case W11(3)}
|f^*(x) - f_{B(x,r)}| \leq C r \alpha,
\end{equation}
where \(C\) is a constant depending only on \(n\).
Indeed, since \(x \notin \mathcal{N}\), we have \(\lim_{r \to 0^+} f_{B(x,r)} = f^*(x) \in \mathbb{R}\). For \(r \in (0, \infty)\) and \(l \in \mathbb{N}\), we apply~\eqref{eq:estimate for differences between averages, in case W11(1)} with \(\frac{r}{2^l}\) in place of \(r\) and \(\frac{r}{2^{l+1}}\) in place of \(s\), to obtain
\begin{equation}
\label{eq:estimate for |f^*(x)-f_B(x,r)|,in case W11}
|f^*(x) - f_{B(x,r)}| = \left| \sum_{l=0}^\infty \left[f_{B(x,r/2^{l+1})} - f_{B(x,r/2^l)} \right] \right|
\leq \sum_{l=0}^\infty \left| f_{B(x,r/2^{l+1})} - f_{B(x,r/2^l)} \right|
\leq C r \alpha.
\end{equation}

Now, we establish the following estimate for the difference between averages of \(f\) over balls with the same radius but different centers, whose distance equals the radius. More precisely, we claim: if \(x, y \in A_\alpha\) with \(x \neq y\), and we set \(r := |x - y|\), then
\begin{equation}
\label{eq:estimate for the averages (f)_x,r and (f)_y,r,in case W11(3)}
|f_{B(x,r)} - f_{B(y,r)}| \leq C r \alpha,
\end{equation}
where \(C\) is a constant depending only on \(n\). Indeed, observe that \(B\left(\frac{x + y}{2}, \frac{r}{2}\right) \subset B(x,r) \cap B(y,r)\), so by the invariance under isometries and the homogeneity of the Lebesgue measure, we have
\begin{equation}
\label{eq: usage of invariance under isometries of Lebesgue measure,in case W11}
\frac{1}{|B(x,r) \cap B(y,r)|} \leq \frac{1}{\left|B\left(\frac{x+y}{2}, \frac{r}{2} \right)\right|}
= \frac{1}{|B(x, \frac{r}{2})|} = \frac{1}{|B(y, \frac{r}{2})|}
= \frac{2^n}{|B(x,r)|} = \frac{2^n}{|B(y,r)|}.
\end{equation}
By the triangle inequality, for almost every \(z \in B(x,r) \cap B(y,r)\), we have
\begin{equation}
\label{eq:triangle inequality for averages,in case W11}
|f_{B(x,r)} - f_{B(y,r)}| \leq |f_{B(x,r)} - f(z)| + |f(z) - f_{B(y,r)}|.
\end{equation}
Taking the average of both sides of~\eqref{eq:triangle inequality for averages,in case W11} over \(B(x,r) \cap B(y,r)\) with respect to \(dz\), and using~\eqref{eq: usage of invariance under isometries of Lebesgue measure,in case W11}, Poincar\'e inequality, and the assumption that $x,y\in A_\alpha$, we obtain
\begin{multline}
\label{eq:estimate for the averages (f)_x,r and (f)_y,r,in case W11}
|f_{B(x,r)} - f_{B(y,r)}|
\leq \fint_{B(x,r) \cap B(y,r)} \left( |f_{B(x,r)} - f(z)| + |f(z) - f_{B(y,r)}| \right) \, dz \\
\leq \frac{2^n}{|B(x,r)|} \int_{B(x,r) \cap B(y,r)} |f_{B(x,r)} - f(z)| \, dz
+ \frac{2^n}{|B(y,r)|} \int_{B(x,r) \cap B(y,r)} |f_{B(y,r)} - f(z)| \, dz 
\\
\leq 2^n \left( \fint_{B(x,r)} |f_{B(x,r)} - f(z)| \, dz
+ \fint_{B(y,r)} |f(z) - f_{B(y,r)}| \, dz \right)
\\
\leq 2^nC(n)r\left( \fint_{B(x,r)} |\nabla f(z)| \, dz
+ \fint_{B(y,r)} |\nabla f(z)| \, dz \right) 
\leq C r \alpha.
\end{multline}
This completes the proof of~\eqref{eq:estimate for the averages (f)_x,r and (f)_y,r,in case W11(3)}.
Combining the estimates
\eqref{eq:estimate for |f^*(x)-f_B(x,r)|,in case W11(3)},
\eqref{eq:estimate for the averages (f)_x,r and (f)_y,r,in case W11(3)},
and using the triangle inequality, we deduce that for every \(x, y \in A_\alpha \setminus \mathcal{N}\) with $x\neq y$ and $r:=|x-y|$,
\begin{equation}
\label{eq:f^* is Lipschitz,in case W11}
|f^*(x) - f^*(y)| \leq |f^*(x) - f_{B(x,r)}| + |f_{B(x,r)} - f_{B(y,r)}| + |f_{B(y,r)} - f^*(y)|
\leq C \alpha |x - y|.
\end{equation}
This completes the proof of the Lipschitz property of \(f^*\) on \(A_\alpha \setminus \mathcal{N}\).
\end{proof}
Now we can prove Theorem \ref{thm:LipCap}.
\begin{proof}[\textbf{Proof of Theorem \ref{thm:LipCap}}]

(1) Assume first 
$f\in W^{k,p}(\mathbb R^n).$
Then for every $\alpha\in (0,\infty)$
we define the set 
\begin{equation}
\label{eq:definition on Aalpha}
A_\alpha=\Set{x\in \mathbb R^n} [M|\nabla f|(x)\leq \alpha].
\end{equation}
By Lemma \ref{lem:level sets of the H.L function are open}, $A_{\alpha}$ is a closed set. Furthermore, by the definition of the sets $A_{\alpha}$, it follows that $A_{\alpha_1}\subset A_{\alpha_2}$ for $\alpha_2\geq\alpha_1>0$.

By Theorem \ref{thm:Lebesgue points of Sobolev functins}, $R_{k-1,p}(E)=0$, where $E$ is the set defined in \eqref{eq:definition of the non Lebesgue points of Ries potential of high derivatives} for $k-1$ in place of $k$, meaning that
\begin{equation}
\label{eq:definition of the non Lebesgue points of Ries potential of high derivatives1}
E:=\Set{x\in\mathbb{R}^n}[I_{k-1}*g(x)=\infty],\quad g(y) := \sum_{|\alpha| = k-1} |D^\alpha f(y)|.
\end{equation}
By Lemma \ref{lem:Lipschitz approximations}, for every $\alpha\in(0,\infty)$, the precise representative $f^*$ is Lipschitz continuous in $A_\alpha\setminus \mathcal{N}$, where the set $\mathcal{N}$ is defined in Lemma \ref{lem:Lipschitz approximations}. Note that $\mathcal{N}\subset E$ because if $x\notin E$, then the limit $\lim_{r\to 0^+}f_{B(x,r)}$ exists and belongs to $\R$, hence $x\notin \mathcal{N}$. Therefore, we get $R_{k-1,p}(\mathcal{N})=0$.

In particular, for every integer
$l\geq1$
the function
$f^*$
is Lipschitz continuous on 
$A_l\setminus \mathcal{N},R_{k-1,p}(\mathcal{N})=0$. By Theorem \ref{thm:weak type inequality for capacities} and Remark \ref{rem:chebyshev inequality for Riesz capacity is valid for vector valued Sobolev mappings}, we get
\begin{equation}
\label{eq:estimate for the capacity by the derivative}
R_{k-1,p}(\mathbb R^n\setminus A_l)\leq \frac{C(n,p,k)}{l^p}\|\nabla f\|_{W^{k-1,p}(\mathbb R^n)}^p.
\end{equation}
Therefore, from
$\eqref{eq:estimate for the capacity by the derivative}$, we conclude that $R_{k-1,p}\left(\mathbb R^n\setminus\bigcup\limits_{l=1}^\infty A_l\right)=0$.
It completes the proof in the case $f\in W^{k,p}(\mathbb R^n)$.

(2) Assume that $f \in W^{k,p}_{\loc}(\Omega)$. Let $\{\Omega_j\}_{j=1}^\infty$ be a nested sequence of open subsets of $\Omega$ such that
\[
\Omega_j \subset \Omega_{j+1} \subset \Omega, \quad \overline{\Omega}_j \subset \Omega_{j+1}, \quad \text{and} \quad \bigcup_{j=1}^\infty \Omega_j = \Omega,
\]
where each $\overline{\Omega}_j$ is compact.

For each $j \in \mathbb{N}$, choose a function $\zeta_j \in C_c^\infty(\Omega_{j+1})$ such that $\zeta_j \equiv 1$ in a neighborhood of $\overline{\Omega}_j$. Define
\[
f_j := f \zeta_j.
\]
Then $f_j \in W^{k,p}(\Omega_{j+1})$ and $\supp(f_j) \subset \Omega_{j+1}$. Extend $f_j$ to $\mathbb{R}^n$ by zero and denote the extension again by $f_j$. It follows that $f_j \in W^{k,p}(\mathbb{R}^n)$.

By the previous case, there exists a non-decreasing sequence of closed sets $\{A_l^j\}_{l=1}^\infty$ (defined as in \eqref{eq:definition on Aalpha}) such that:
 $f_j^*|_{A_l^j}$ is Lipschitz continuous on $A_l^j$ except for a set of $R_{k-1,p}$-capacity zero,
 and $R_{k-1,p}\left( \mathbb{R}^n \setminus \bigcup_{l=1}^\infty A_l^j \right) = 0$.

Let us now choose a sequence $\{\alpha_j\}_{j=1}^\infty \subset \mathbb{N}$ such that
\[
\alpha_j^p \geq 2^j \| \nabla f_j \|_{W^{k-1,p}(\mathbb{R}^n)}^p.
\]
Then, by Theorem~\ref{thm:weak type inequality for capacities}, we obtain
\begin{equation} \label{eq:estimate for the Riesz capacity of the complement of A}
R_{k-1,p}(\mathbb{R}^n \setminus A_{\alpha_j}^j)
\leq \frac{C(n,p,k)}{\alpha_j^p} \| \nabla f_j \|_{W^{k-1,p}(\mathbb{R}^n)}^p 
\leq \frac{C(n,p,k)}{2^j}.
\end{equation}

Define a sequence of sets $\{B_l\}_{l=1}^\infty$ by
\[
B_l := \bigcap_{j=l}^\infty A_{\alpha_j}^j.
\]
Then, the complement of their union satisfies
\begin{equation} \label{eq:representaion of the complement of the union of Bl}
\mathbb{R}^n \setminus \bigcup_{l=1}^\infty B_l 
= \bigcap_{l=1}^\infty \bigcup_{j=l}^\infty \left( \mathbb{R}^n \setminus A_{\alpha_j}^j \right).
\end{equation}
Thus, for any $i \in \mathbb{N}$, using \eqref{eq:representaion of the complement of the union of Bl} and \eqref{eq:estimate for the Riesz capacity of the complement of A}, we have
\begin{equation} \label{eq:complement of union of Bl has capacity zero1}
R_{k-1,p}\left( \mathbb{R}^n \setminus \bigcup_{l=1}^\infty B_l \right) 
\leq R_{k-1,p}\left( \bigcup_{j=i}^\infty \left( \mathbb{R}^n \setminus A_{\alpha_j}^j \right) \right)
\leq \sum_{j=i}^\infty R_{k-1,p}\left( \mathbb{R}^n \setminus A_{\alpha_j}^j \right)
\leq C(n,p,k) \sum_{j=i}^\infty \frac{1}{2^j}.
\end{equation}
Since this holds for all $i \in \mathbb{N}$, we conclude that
\begin{equation} \label{eq:complement of union of Bl has capacity zero}
R_{k-1,p}\left( \mathbb{R}^n \setminus \bigcup_{l=1}^\infty B_l \right) = 0.
\end{equation}

Now define
\[
C_l := B_l \cap \overline{\Omega}_l.
\]
The sequence $\{C_l\}_{l=1}^\infty$ is non-decreasing, since both $\{B_l\}$ and $\{\overline{\Omega}_l\}$ are non-decreasing. Each $C_l$ is closed as the intersection of closed sets. Moreover,
\[
\bigcup_{l=1}^\infty C_l = \left( \bigcup_{l=1}^\infty B_l \right) \cap \Omega.
\]
Therefore,
\begin{equation} 
\label{eq:complement of union of Cl lies in complement of Bl}
\Omega \setminus \bigcup_{l=1}^\infty C_l
=\Omega \setminus \left( \bigcup_{l=1}^\infty B_l \right)
\subset \mathbb{R}^n \setminus \bigcup_{l=1}^\infty B_l.
\end{equation}
Hence, combining \eqref{eq:complement of union of Cl lies in complement of Bl} and \eqref{eq:complement of union of Bl has capacity zero}, we obtain
\[
R_{k-1,p}\left( \Omega \setminus \bigcup_{l=1}^\infty C_l \right) = 0.
\]

Finally, note that for each $l \in \mathbb{N}$, we have
\[
C_l \subset A_{\alpha_l}^l \cap \overline{\Omega}_l,
\]
and that $f_l^*|_{A_{\alpha_l}^l \cap \overline{\Omega}_l} = f^*|_{A_{\alpha_l}^l \cap \overline{\Omega}_l}$, since $\zeta_l \equiv 1$ on a neighborhood of $\overline{\Omega}_l$. Since $f_l^*$ is Lipschitz continuous on $A_{\alpha_l}^l \cap \overline{\Omega}_l$ outside a set of $R_{k-1,p}$-capacity zero, it follows that $f^*|_{C_l}$ is Lipschitz continuous outside a set of $R_{k-1,p}$-capacity zero.

This completes the proof of the theorem. 
\end{proof}

\begin{remark}
Note that for every $k\in \mathbb{N}$ and $1\leq p\leq \infty$, and for every $f\in W^{k,p}(\mathbb{R}^n)$, Theorem \ref{thm:LipCap} follows from Lemma \ref{lem:Lipschitz approximations} if we replace the Riesz capacity with Lebesgue measure. Indeed, by Lemma \ref{lem:Lipschitz approximations}, the precise representative $f^*$ is Lipschitz continuous on $A_\alpha \setminus \mathcal{N}$, where $A_\alpha$ and $\mathcal{N}$ are defined as in Lemma \ref{lem:Lipschitz approximations}. By the Lebesgue differentiation theorem, $|\mathcal{N}| = 0$. Furthermore, by the Hardy–Littlewood maximal inequality, there exists a constant $C = C(n,p)$ depending only on $n$ and $p$ such that  
\begin{equation}
\label{eq:Lebesgue measure of the complement of A}
\left|\mathbb{R}^n \setminus A_\alpha\right|= \left|\left\{x\in\mathbb{R}^n\,:\,M|\nabla f|(x)>\alpha\right\}\right| \leq \frac{C}{\alpha^p} \|\nabla f\|^p_{L^p(\mathbb{R}^n)}.
\end{equation}
Choosing any sequence of positive numbers $\alpha_j$ that converges to infinity as $j \to \infty$, we obtain a monotonically non-decreasing sequence $\{A_{\alpha_j}\}$ of closed sets such that $f^*$ is Lipschitz continuous on $A_{\alpha_j}$ outside a set of Lebesgue measure zero, and we get from \eqref{eq:Lebesgue measure of the complement of A} 
\begin{equation}
\left|\mathbb{R}^n \setminus \bigcup_{j=1}^\infty A_{\alpha_j}\right| = 0.
\end{equation}
The general case $f \in W^{k,p}_{\loc}(\Omega)$, where $\Omega\subset\mathbb{R}^n$ is an open set, is proved similarly to the proof of part (2) of Theorem \ref{thm:LipCap}.  
\end{remark}

\begin{proof}[\textbf{Proof of Theorem \ref{thm:area formula for Riesz capacity}}]
Theorem~\ref{thm:area formula for Riesz capacity} follows by combining Theorem~\ref{thm:LipCap} and Theorem~\ref{thm:general area formula}. 

From Theorem~\ref{thm:LipCap}, we obtain a sequence of closed (and in particular Lebesgue measurable) sets $C_l$ such that $C_l \subset C_{l+1}$ and $\varphi^*$ is Lipschitz continuous on $C_l$ $R_{k-1,p}$-almost everywhere, for every $l \in \mathbb{N}$. Moreover,
\[
R_{k-1,p}\left(\Omega \setminus \bigcup_{l=1}^\infty C_l\right) = 0.
\]
For each $l \in \mathbb{N}$, let $\Theta_l \subset C_l$ be the $R_{k-1,p}$-null set such that $\varphi^*$ is Lipschitz continuous on $C_l \setminus \Theta_l$. Set 
\[
\Theta := \bigcup_{l=1}^\infty \Theta_l, \quad \text{and} \quad E_l := C_l \setminus \Theta.
\]
Then the sequence of sets $E_l$ satisfies the assumptions of Theorem~\ref{thm:general area formula}. 

Let $f : \Omega \to \mathbb{R}$ be any non-negative Lebesgue measurable function, and define 
\[
S := \Omega \setminus \bigcup_{l=1}^\infty E_l.
\]
By Theorem~\ref{thm:general area formula}, we obtain:
\begin{equation}
\label{eq:area formula with S}
\int_{\Omega \setminus S} f(x) |J(x, \varphi)| \, dx = \int_{\mathbb{R}^n} \left( \sum_{x \in \varphi^{-1}\{y\} \cap (\Omega \setminus S)} f(x) \right) \, dy.
\end{equation}

By Remark~\ref{rem:Riesz capacity is an outer measure}, we know that $R_{k-1,p}(S) = 0$ implies $|S| = 0$. Therefore, the set $S$ may be omitted from the left-hand side of~\eqref{eq:area formula with S}, and we conclude:
\begin{equation}
\label{eq:area formula without S}
\int_{\Omega} f(x) |J(x, \varphi)| \, dx = \int_{\mathbb{R}^n} \left( \sum_{x \in \varphi^{-1}\{y\} \cap (\Omega \setminus S)} f(x) \right) \, dy.
\end{equation}
\end{proof}

\section*{Funding}
This work was supported by the Israel Science Foundation (grant No.\ 569/21).

\begin{acknowledgment}
The author would like to express his gratitude to Alexander Ukhlov, who suggested the topic and provided helpful suggestions during the preparation of this article.
\end{acknowledgment}

\textbf{Data Availability}: This article does not involve any datasets.

\textbf{There is no conflict of interest}

\vskip 0.3cm

Department of Mathematics, Ben-Gurion University of the Negev, P.O.Box 653, Beer Sheva, 8410501, Israel

\emph{E-mail address:} \email{pazhash@post.bgu.ac.il}

\end{document}